\documentclass[10pt]{article}
\usepackage{amsmath}
\usepackage{amsthm}
\usepackage{amssymb}
\usepackage{amscd}
\usepackage{amsfonts}
\usepackage{amssymb}
\usepackage{setspace}
\setdisplayskipstretch{0.8}

 \newtheorem{theorem}{Theorem}[section]
 
 \newtheorem{lemma}[theorem]{Lemma}
 
 \newtheorem{proposition}[theorem]{Proposition}

\theoremstyle{definition}
\newtheorem{example}{Example}

 \def\beqa{\begin{eqnarray*}}
\def\eeqa{\end{eqnarray*}}

\newcommand{\D}{\Delta}
\newcommand{\mflr}{\left\lfloor \frac{m}{2} \right\rfloor}
\newcommand{\mcl}{\left\lceil \frac{m}{2} \right\rceil}

\begin{document}

\title{On a modification of a problem of Bialostocki, Erd\H{o}s, and Lefmann}

\author{Andrew Schultz\thanks{This research was
supported in part by funding from NSF grant DMS0097317}}
\maketitle
\begin{abstract}
For positive integers $m$ and $r$, one can easily show there exist integers $N$ such that for every map $\D:\{1,2,\cdots,N\} \to \{1,2,\cdots,r\}$ there exist $2m$ integers $$x_1
< \cdots < x_m < y_1 < \cdots < y_m$$ which satisfy:
\begin{enumerate} \item $\hspace{81 pt}\D(x_1) = \cdots = \D(x_m)$, \item
$\hspace{81 pt} \D(y_1) = \cdots = \D(y_m)$, \item $\hspace{81 pt}
2(x_m-x_1) \leq y_m-x_1$.
\end{enumerate}

In this paper we investigate the minimal such integer, which we call $g(m,r)$.  We prove
that $g(m,2) = 5(m-1)+1$ for $m \geq 2$, that $g(m,3) =
7(m-1)+1+\mcl$ for $m \geq 4$, and that $g(m,4)=10(m-1)+1$ for $m
\geq 3$.  Furthermore, we consider $g(m,r)$ for general $r$. Along with results that bound $g(m,r)$, we compute $g(m,r)$ exactly for the following infinite families of $r$: $$\left\{f_{2n+3}\right\}, \left\{2f_{2n+3}\right\}, \left\{18f_{2n}-7f_{2n-2}\right\}, \mbox{ and } \left\{23f_{2n}-9f_{2n-2}\right\}, $$ where here $f_i$ is the $i$th Fibonacci number defined by $f_0 = 0$ and $f_1=1$.
\end{abstract}
\onehalfspacing
\section{Introduction}

Ramsey type problems regarding colorings of the natural numbers
are concerned with finding the minimum number $N(r)$, if it exists, for which
every coloring of the integers in $[1,N]$ by $r$ colors contains
some given monochromatic configuration.  Traditionally, these
configurations are solutions to systems of linear equations.  The
general theory developed by Rado in \cite{16} gave rise to the
determination of $N(r)$ for certain monochromatic configurations,
such as Schur numbers and Van der Waerden numbers \cite{Gr}. Other
exact results of a similar kind were determined in \cite{13},
\cite{14}, \cite{17}, \cite{18}, and \cite{19}. The difficulty in
computing such numbers led to the consideration of inequalities
instead of equations.  In particular, arithmetic progressions
prompted Brown, Erd\H{o}s, and Freedman to define the notion of
ascending waves.  These and similar structures have been
investigated in \cite{1}, \cite{7}, and \cite{15}.

Along similar lines, Bialostocki, Erd\H{o}s, and Lefmann
considered in \cite{3} the following problem concerning
monochromatic sets of nondecreasing diameter.  For two positive integers $m$ and $r$, determine the minimum integer,
$f(m,r)$, such that for every map $\D:\{1, \cdots, f(m,r) \to \{1,\cdots,r\}$ there exist $2m$
integers $$x_1 < \cdots < x_m < y_1 < \cdots < y_m$$ which satisfy
the following conditions:
\begin{enumerate} \item $\hspace{97 pt}\D(x_1) = \cdots =
\D(x_m)$, \item $\hspace{97 pt}\D(y_1) = \cdots = \D(y_m)$, \item
$\hspace{97 pt}x_m-x_1 \leq y_m-y_1$.
\end{enumerate}  They showed $f(m,2)=5m-3$ and $f(m,3)=9m-7$.
 Recently, Grynkiewicz proved $f(m,4)=12m-9$ in \cite{12}.
 Bollob\'{a}s, Erd\H{o}s, and Jin investigated in \cite{5} a
 closely related function, $f^*(2,r)$, where strict inequality is
 required in 3 above.  They determined $f^*(2,r)$
 for $r=2^k$.

In this paper we replace condition 3 by
\begin{enumerate} \item[3$'$.] $\hspace{97 pt}2(x_m-x_1) \leq
y_m-x_1$\end{enumerate} and denote the corresponding function by
$g(m,r)$.  Notice that this is a relaxation of 3, since adding 
$x_m-x_1$ to each side of the inequality $x_m-x_1 \leq y_m-y_1$ 
yields $2(x_m-x_1)$ on the left and $y_m-y_1+x_m-x_1 < y_m-x_1$ 
on the right.  It is not hard to see $g(1,r)=2$ for any $r \in \mathbb{N}$, and 
as such we assume throughout the sequel that $m \geq 2$.

The paper is organized as follows.  In Section 2 we introduce some
basic terms and develop a useful lemma that simplifies the
construction of lower bounds.  In Section 3 we determine $g(m,r)$
for $r \in \{2,3,4\}$.  The main theorems appear in Section 4,
where we develop tools that allow for either a bound or
determination of $g(m,r)$ based upon the value of $g(m,j)$ for
$j<r$.  We conclude with some conjectures that arose from
studying $g(m,r)$ using a computer program based on the theorems
in Section 4.

\section{Preliminaries}

If $S$ is a nonempty set of integers and $\D:S \rightarrow R$ is a
mapping where $|R|=r$, then $\D$ is called an $r$-coloring of $S$.  
For $T \subseteq S$, we write $\Delta(T) = \{\Delta(t): t \in T\}$. 
We say $T$ is monochromatic if $|\D(T)|=1$.

Throughout this paper an $m$-set, denoted $Z=(z_1,\ldots,z_m)$, is
a sequence of $m$ distinct positive integers such that
$z_1<\cdots<z_m$.  For a pair of $m$-sets $X$ and $Y$, we write $X
\prec Y$ if $x_m<y_1$.  Suppose $X \prec Y$; we define $Y$ to be
$X$-admissible if $2(x_m-x_1)>y_m-x_1$.  Furthermore, let $\D$ be
an $r$-coloring of a nonempty set $S$; we say $\D$ is an
$L(r)$-coloring of $S$ if for every pair of monochromatic $m$-sets
$X,Y \subset S$, either $X \nprec Y$ or $Y$ is $X$-admissible.
That is, a coloring $\D$ is an $L(r)$-coloring provided there are
no two monochromatic $m$-sets $X,Y \subset S$ such that $X \prec
Y$ and conditions 1,2, and 3$'$ above are satisfied.

For an $n$-set $X=(x_1,\ldots,x_n)$ we use the following notation:

\noindent(i) $int_i(X)=x_i$ for $i \leq n$;
\newline
(ii)  $first_k(X)=\{x_1,\ldots,x_{\min\{k,n\}}\}$; and
\newline
(iii)  $last_k(X)=\{x_{\max\{1,n-k+1\}},\ldots,x_n\}$.

For two integers $a$ and $b$ we use $[a,b]$ to
denote the set of all integers $i$ such that $a \leq i$ and 
$i \leq b$, and refer to it as an interval.  Note that if $a>b$ 
then $[a,b]=\varnothing$.  Furthermore, for positive
integers $r$,$m$, and $s$, where $s \geq 2r(m-1)+1$, define the
disjoint intervals $I_1$, $I_2$, and $I_3$ to be \begin{enumerate}
\item $I_1 = [1,r(m-1)+1]$, \item $I_2 = [r(m-1)+2,2r(m-1)]$, and
\item $I_3=[2r(m-1)+1,s]$.\end{enumerate}  Here we have used $m \geq 2$ to assume $I_1 \cap I_3 =\varnothing$.

Since $|I_1|=r(m-1)+1$ one sees that for an arbitrary $r$-coloring $\D$ there must be some monochromatic $m$-set $X \subseteq I_1$.
The following proposition is immediate.

\begin{proposition}\label{prop:first}  Let $s \geq 2r(m-1)+1$ be a positive integer, and let $\D: [1,s] \rightarrow [1,r]$
be a coloring.  If there exists a monochromatic $m$-set $Y \subset
I_2 \cup I_3$ with $y_m \in I_3$ then $\D$ is not an
$L(r)$-coloring.
\end{proposition}

The following lemma simplifies the construction of
$L(r)$-colorings by inducing an $L(r)$-coloring of $I_1 \cup I_2$
from an $L(r)$-coloring of $I_2$.

\begin{lemma}\label{lem:colorextension}  Let $\D: I_2 \rightarrow [1,r]$ be an
$L(r)$-coloring.  Then there exists an $L(r)$-coloring $\D_e$ of
$I_1 \cup I_2$ which is an extension of $\D$.  Further, $\Delta_e$ satisfies 
$\Delta_e(1) = \Delta_e(r(m-1)+1)$ and $|\Delta_e^{-1}(t)\cap[1,r(m-1)]|=m-1$ for all $t$.
\end{lemma}

\begin{proof}
Since $|I_2|=r(m-1)-1$, it follows that there is a color $c$ such
that $|\D^{-1}(c)\cap I_2|<m-1$.

We define $\D_e: I_1 \cup I_2 \rightarrow [1,r]$ in two steps.
First, we induce a coloring on $I_2$ and part of $I_1$ as
described below:\[ \D_e(x) =
   \left\{ \begin{array}{ll}
       c, &\mbox{if }x=1 \mbox{ or }x=r(m-1)+1 \\
       \D\left(x+r(m-1)\right), &\mbox{if } x+r(m-1) \in \cup_{t=1}^r first_{m-1} (\D^{-1}(t)\cap I_2)\\
       \D(x), &\mbox{if } x \in I_2 \end{array} \right. \]

Second, we color the remaining integers of $I_1$ recursively as
follows:  suppose $x\in I_1$ and that $\D_e \big{|}_{[1,x-1]}$ is
defined while $\D_e(x)$ is not; then $\D_e(x)=i$, where
$i=\min{[1,r]}$ such that $|\D_e^{-1}(i) \cap [1,r(m-1)]|<m-1$.
From the definition of $\D_e$ it is easy to verify that
$|\D_e^{-1}(t)\cap [1,r(m-1)]|=m-1$ for every $t \in [1,r]$.

It is left to show that $\D_e$ is an $L(r)$-coloring of $I_1 \cup
I_2$.  Let $X, Y \subset I_1 \cup I_2$ be monochromatic $m$-sets
with $X \prec Y$.  If $x_1 \in I_2$ then $Y$ is $X$-admissible
since $\D_e \big{|}_{I_2} = \D$, which is an $L(r)$-coloring by
assumption.  Hence we may assume $x_1 \in I_1$.

\textbf{Case 1}.  Suppose $\D_e(x_1) =t \neq c$.  Then since
$|\D_e^{-1}(t) \cap I_1|=m-1$ for each $t \neq c$, it follows that
$x_1=int_i(\D_e^{-1}(t)\cap I_1)$ for some $i \in [1,m-1]$. Hence,
since $\D_e(r(m-1)+1)=c \neq t$, it follows that
$|\D_e^{-1}(t)\cap [x_1,r(m-1)+1]|=m-i$, and thus $x_m \geq
int_i(\D_e^{-1}(t)\cap I_2)$ and $|\D^{-1}(t)\cap I_2|\geq i$. Remembering $i \leq m-1$, it therefore follows from the
definition of $\D_e$ that
$$x_1=int_i(\D_e^{-1}(t)\cap I_1) = int_i(\D^{-1}(t)\cap I_2) -r(m-1) \leq x_m -r(m-1),$$ so that $x_m-x_1 \geq r(m-1)$.
Hence, since $Y \subset I_1 \cup I_2$, $$2(x_m-x_1) \geq 2r(m-1)
\geq y_m > y_m-x_1,$$ and $Y$ is $X$-admissible.

\textbf{Case 2}.  Suppose $\D_e(x_1)=c$.  The argument above holds
except in the case that $x_1=1$.  In this case, we have $x_m \geq
r(m-1)+1=r(m-1)+x_1$, and $x_m-x_1 \geq r(m-1)$ as before.
\end{proof}

In conjunction with Proposition \ref{prop:first}, Lemma
\ref{lem:colorextension} shows there exists an $L(r)$-coloring
on $I_1 \cup I_2 \cup I_3$ provided the existence of a coloring
$\D: I_2 \cup I_3 \rightarrow [1,r]$ which is an $L(r)$-coloring
on $I_2$ such that $|\D^{-1}(c) \cap (I_2 \cup I_3)| \leq m-1$
for every $c \in \D(I_3)$.  Henceforth, we shall let the existence
of $\D: I_2 \cup I_3 \rightarrow [1,r]$ which satisfies these
conditions suffice to show the existence of an $L(r)$-coloring
$\D_e: I_1 \cup I_2 \cup I_3 \rightarrow [1,r]$ without explicit
construction.

\section{The function $g(m,r)$ for $r \in \{2,3,4\}$}

We first evaluate the function $g(m,r)$ for small values of $r$
and appropriate values of $m$.  The case when $r=2$ is trivial.

\begin{theorem} \label{thm:g(m,2)}Let $m \geq 2$ be an integer.  Then, $g(m,2) =
5m-4$.\end{theorem}
\begin{proof}
The coloring $\D: [2m,5m-5] \rightarrow [1,2]$ given by
$$1^{2m-3}2^{m-1}$$ shows that $g(m,2) \geq 5m-4$.

Next we show that $g(m,2) \leq 5m-4$.  Let $\D: [1,5m-4]
\rightarrow [1,2]$ be an arbitrary $2$-coloring, and let
$P=[3m-2,5m-4]$.  Since $|P|=2m-1$ there exists some monochromatic
$m$-set $Y \subset P$.  Furthermore, since $|P \cap I_2|=m-1$, it
follows that $Y \cap I_3 \neq \varnothing$. Applying Proposition
\ref{prop:first} completes the proof.
\end{proof}

In evaluating $g(m,3)$ it will be beneficial to have the following

\begin{lemma}\label{lem:3colorlemma}
Let $m \geq 4$ be an integer, and let $\D: [1,3m-4]
\rightarrow [1,3]$ be a $3$-coloring.  If $|\D^{-1}(c)| \geq 3m-
\big{\lceil} \frac{m}{2} \big{\rceil}-2$ for some $c \in [1,3]$,
then $\D$ is not an $L(3)$-coloring.
\end{lemma}
\begin{proof}  Let $I=[1,3m-4]$ and $t=|[1,int_1(\D^{-1}(c)\cap I)-1]|$. Further, let $s=3m-4-|\D^{-1}(c)| \leq \lceil
\frac{m}{2} \rceil-2$, the number of integers in the interval $I$
not colored by $c$. Finally, let $w = |[int_1(\D^{-1}(c)\cap I),
int_m(\D^{-1}(c)\cap I)]|-m$.  It will be important later to note
that
\begin{equation}\label{eqn:1} w+s \leq 2s \leq 2(\left\lceil
\frac{m}{2}\right\rceil-2) \leq m-3.
\end{equation}

Let $X = first_m(\D^{-1}(c)\cap I)$ (note that since $\D^{-1}(c)
\geq 3m-\lceil \frac{m}{2} \rceil -2 > m$, $X$ is in fact an
$m$-set). By construction we have $x_1 = t+1$ and $x_m = t+w+m$,
so that $x_m-x_1 = m+w-1$.  Hence, if there is a monochromatic
$m$-set $Y$ with $y_m \geq x_1+2(m-1+w) = 2m-1+t+2w$ and $X \prec
Y$, then $Y$ is not $X$-admissible and the proof is complete.  We
show that $$Y = last_m(\D^{-1}(c)\cap I)$$ satisfies these
conditions.

First, note that $|\D^{-1}(c)| \geq 3m-\lceil \frac{m}{2} \rceil
-2 \geq 2m$ since $m \geq 4$, from which it follows that $Y$ is
indeed an $m$-set and $X \prec Y$.  We now show
$last_1(\D^{-1}(c)\cap I) \geq 2m-1+t+2w$.  Since there are
exactly $s-(t+w)$ integers $z$ with $\D(z) \neq c$ and $z>x_m$, it
follows that $last_1(\D^{-1}(c)\cap I) \geq
3m-4-\left(s-(t+w)\right)$.  Hence, recalling Equation
\ref{eqn:1}, it follows that
\begin{equation*}\begin{split}last_1(\D^{-1}(c)\cap I) &\geq
3m-4-s+t+w \\ &\geq 3m-4-(m-3-w)+t+w
\\ &=2m-1+t+2w,\end{split}\end{equation*} and the proof is complete.
\end{proof}

\begin{theorem}\label{thm:g(m,3)}Let $m \geq 4$ be an integer.  Then, $g(m,3) = 7m +
\mcl -6$.
\end{theorem}
\begin{proof}
One may verify that the coloring $\D: [3m-1,7m + \mcl
-7] \rightarrow [1,3]$ given by
$$1^{m-\mflr-2}2^{\mflr-1}1^{2m-1}2^{\mcl}3^{m-1}$$ shows $g(m,3)
\geq 7m + \mcl -6$.

Next we show that $g(m,3) \leq 7m + \mcl -6$.  Let $\D: [1,7m +
\mcl -6] \rightarrow [1,3]$ be an arbitrary $3$-coloring.  Since
$|I_2| = 3(m-1)-1$ it follows there exists some $c \in [1,3]$ such
that $|\D^{-1}(c) \cap I_2| \geq m-1$.  If $\D^{-1}(c)\cap I_3
\neq \varnothing$, then the proof is complete.  We may therefore
assume $\D^{-1}(c)\cap I_3 = \varnothing$ and thus $|\D(I_3)|\leq
2$.  Since $|I_3|=m+\mcl$, if $|I_2|-|\D^{-1}(c) \cap I_2| \geq
\mflr-1$ then it follows from the pigeonhole principle that some
monochromatic $m$-set $Y \subset I_2 \cup I_3$ exists with $Y \cap
I_3 \neq \varnothing$.  In this case, an application of
Proposition \ref{prop:first} completes the proof.

Finally, we are left to assume that $|I_2|-|\D^{-1}(c) \cap I_2| <
\mflr-1$, so that $|\D^{-1}(c) \cap I_2| \geq 3m-\mcl-3$.  Translating
$I_2$ to the interval $[1,3m-4]$ and applying Lemma
\ref{lem:3colorlemma} completes the proof.
\end{proof}

\begin{lemma}\label{lem:4colorlemma}Let $m \geq 3$ be an integer, and let $\D: [1,4m-5]
\rightarrow [1,4]$ be a $4$-coloring.  If $|\D^{-1}(c)| \geq 3m-3$
for some $c \in [1,4]$, then $\D$ is not an $L(4)$-coloring.
\end{lemma}

\begin{proof}The proof of Lemma \ref{lem:4colorlemma} is similar
to that of Lemma \ref{lem:3colorlemma}, and we omit it.\end{proof}

\begin{theorem}\label{thm:g(m,4)}Let $m \geq 3$ be an integer.  Then,
$g(m,4)=10m-9$.\end{theorem}
\begin{proof}One may verify that the coloring $\D: [4m-2,10m-10]
\rightarrow [1,4]$ given by
$$1^{m-3}2^{m-1}1^{2m-1}3^{m-1}4^{m-1}$$ shows $g(m,4) \geq
10m-9$.

Next we show that $g(m,4) \leq 10m-9$.  Let $\D: [1,10m-9]
\rightarrow [1,4]$ be an arbitrary $4$-coloring.  Since
$|I_2|=4(m-1)-1$, it follows that there exists $c \in [1,4]$ such
that $|\D^{-1}(c) \cap I_2| \geq m-1$.  If $\D^{-1}(c) \cap I_3
\neq \varnothing$ then the proof is complete.  Otherwise we have
$\D^{-1}(c) \cap I_3 = \varnothing$, and so $|\D(I_3)|\leq 3$.
Since $|I_3|=2m-1$, if $|I_2|-|\D^{-1}(c) \cap I_2| \geq m-1$ then
it follows that some monochromatic $m$-set $Y \subset I_2 \cup
I_3$ exists with $Y \cap I_3 \neq \varnothing$.  In this case, an
application of Proposition \ref{prop:first} completes the proof.

Finally, we are left to assume that $|I_2|-|\D^{-1}(c) \cap I_2| <
m-1$, so that $|\D^{-1}(c) \cap I_2| \geq 3m-3$.  Translating $I_2$
to the interval $[1,4m-5]$ and applying Lemma
\ref{lem:4colorlemma} completes the proof.
\end{proof}

\section{Recursion in evaluating $g(m,r)$ when $r \geq
5$}

Though the techniques used in the previous section may be
duplicated in an attempt to solve $g(m,r)$ for $r > 4$, the
limitations of such an approach are easily seen.  In this section
we instead focus our attention on a more general argument which
will allow us to solve $g(m,r)$ for certain infinite families of
integers.

Developing this technique will require that we know certain
properties of $L(r)$-colorings.  The following two lemmas give
some information concerning the structure of $L(r)$-colorings on
the interval $[1,g(m,r)-k]$.

\begin{lemma}\label{lem:nottoomanytostartwith}
Let $m,r$ be positive integers, and let $\D: [1,g(m,r)-1]
\rightarrow [1,r]$ be an $r$-coloring.  If $|\D^{-1}(c)\cap
[1,r(m-1)]| \geq m$ for some $c \in [1,r]$, then $\D$ is not an
$L(r)$-coloring.\end{lemma}
\begin{proof}Suppose for contradiction's sake that $X$ is a monochromatic $m$-set with
$X \subset [1,r(m-1)]$ and that $\D$ is an $L(r)$-coloring of
$[1,g(m,r)-1]$.  Then $\D \big{|}_{[x_m+1,g(m,r)-1]}$ is an
$L(r)$-coloring such that no monochromatic $m$-set $Y$ exists with
$Y \subset [x_m+1,g(m,r)-1]$ and $Y \cap [2x_m-x_1,g(m,r)-1] \neq
\varnothing$.  Applying Lemma \ref{lem:colorextension} and its
subsequent remark allows us to extend $\D
\big{|}_{[x_m+1,g(m,r)-1]}$ to an $L(r)$-coloring $\D_e$ of the
interval $[x_m-r(m-1),g(m,r)-1]$.  Since $x_m \leq r(m-1)$, it
follows that $\D_e \big{|}_{[0,g(m,r)-1]}$ is an $L(r)$-coloring,
which after an appropriate translation contradicts the definition of $g(m,r)$.
\end{proof}

\begin{lemma}\label{lem:needsomeofeachcolor}Let $m,r$ and $k$ be positive integers, and let $\D:
[1,g(m,r)-k] \rightarrow [1,r]$ be an $r$-coloring.  Let
$a=\min\left\{int_m(\D^{-1}(c))\right\}_{c \in [1,r]}$. For each $c \in
[1,r]$, let $A_c(\D)=|\D^{-1}(c) \cap [1,a-1]|$ and
$B_c(\D)=|\D^{-1}(c) \cap [a+1,g(m,r)-k]|$. If
\begin{equation}\label{eq:needsomeofeachcolor}\sum_{c \in [1,r]}\left(A_c(\D)+ \min\{B_c(\D),m-1\}\right) \leq r(2m-2)-k\end{equation} then $\D$ is not an
$L(r)$-coloring.
\end{lemma} \begin{proof}We use induction on $k$.  Suppose $k=1$, and assume for contradiction's sake that $\Delta$ is an $L(r)$-coloring.  By Lemma
\ref{lem:nottoomanytostartwith}, it must be the case that
$a=r(m-1)+1$, so that $[1,a]=I_1$ and $[a+1,g(m,r)-1]=I_2 \cup
I_3$.  Hence we have $$\sum_{c \in [1,r]} A_c(\D)=r(m-1),$$ so by Equation \ref{eq:needsomeofeachcolor} it
must be the case that $|\D^{-1}(c) \cap (I_2 \cup I_3)|<m-1$ for
some $c \in [1,r]$. Induce a coloring $\D_e: [1,g(m,r)]
\rightarrow [1,r]$ defined by
\[ \D_e(x) =
   \left\{ \begin{array}{ll}
       \D(x), &\mbox{ for }x \in [1,g(m,r)-1] \\
       c, &\mbox{ for } x =g(m,r). \end{array} \right. \]

By the definition of $g(m,r)$ there exist $m$-sets $X,Y \subset
[1,g(m,r)]$ with $X \prec Y$ and $y_m-x_1 \geq 2(x_m-x_1)$. Since
$\D$ is an $L(r)$-coloring, it follows that
$y_m=g(m,r)$; furthermore $y_1 \in I_1$ since $|\D^{-1}(c) \cap
(I_2 \cup I_3)| < m-1$.  Therefore, $X \subset [1,r(m-1)]$, a
contradiction.

Assume the result holds for $k$; we show it also holds for $k+1$.
Let $\D:[1,g(m,r)-k-1] \rightarrow [1,r]$ be such that
\begin{equation}\label{eq:inductiveassump}\sum_{c \in [1,r]}A_c(\D)+
\min\{B_c(\D),m-1\}\leq r(2m-2)-k-1.\end{equation}  We consider
two cases.

\textbf{Case 1.}  If $a<r(m-1)+1$, then there must be some $t \in
[1,r]$ such that $|\D^{-1}(t) \cap [1,a]|<m-1$. Induce a coloring
$\D_e: [1,g(m,r)-k] \rightarrow [1,r]$ defined by
\[ \D_e(x) =
   \left\{ \begin{array}{ll}
       t, &\mbox{ for } x =1  \\
       \D(x-1), &\mbox{ for }x \in [2,g(m,r)-k] \end{array} \right. \]
Notice that for $\D_e$ we have
$\min\{int_m(\D^{-1}(c)\cap[1,g(m,r)-k])\}_{c \in [1,r]} = a+1$ so
that
$$\sum_{c \in [1,r]}A_c(\D_e)+
\min\{m-1,B_c(\D_e)\} \leq r(2m-2)-k.$$  Hence, by induction there
exist monochromatic $m$-sets $X,Y$ with $X \prec Y$ and
\begin{equation}\label{eq:thecondition}y_m-x_1 \geq 2(x_m-x_1).\end{equation}  If $\D$ is an $L(r)$-coloring
it follows that $x_1=1$; furthermore $x_m > a+1$ since
$|\D_e^{-1}(t) \cap [1,a+1]| \leq m-1$.  Denoting the
monochromatic $m$-set $first_m (\D_e(a+1)^{-1}\cap[1,a+1])$ by
$Z$, we therefore have $x_1<z_1$ and $x_m > z_m$.  Along with
Equation \ref{eq:thecondition}, this gives us
$$y_m+z_1 > y_m+1 \geq 2x_m > 2z_m,$$ from which it follows that $y_m-z_1 \geq 2(z_m - z_1)$, a contradiction.  Therefore, $\D$ is not an
$L(r)$-coloring.

\textbf{Case 2.} If $a \geq r(m-1)+1$ (and hence $a=r(m-1)+1$), we have
\begin{equation}\label{eq:notmuchofeachcolor}|\D^{-1}(c) \cap [1,r(m-1)]| = m-1 \end{equation} for every $c \in
[1,r]$.  By Equation \ref{eq:inductiveassump}, there
must be some $t \in [1,r]$ such that $|\D^{-1}(t)\cap (I_2 \cup
I_3)|<m-1$.  Induce a coloring $\D_e: [1,g(m,r)-k] \rightarrow
[1,r]$ defined by
\[ \D_e(x) =
   \left\{ \begin{array}{ll}
       \D(x), &\mbox{ for }x \in [1,g(m,r)-k-1] \\
       t, &\mbox{ for } x =g(m,r)-k.  \end{array} \right. \]
It is easily verified for $\D_e$ that
$$\sum_{c \in [1,r]}A_c(\D_e)+
\min\{B_c(\D_e),m-1\} \leq r(2m-2)-k.$$  Hence, by induction there
exist monochromatic $m$-sets $X,Y$ with $X \prec Y$ such that
$y_m-x_1 \geq 2(x_m-x_1)$.  If $\D$ is an $L(r)$-coloring it
follows that $y_m=g(m,r)-k$; furthermore $y_1 \leq
r(m-1)+1$ since $|\D_e^{-1}(t) \cap (I_2 \cup I_3)| \leq m-1$.
Hence, $X \subset [1,r(m-1)]$, a contradiction.\end{proof}

We now develop a recursive technique for evaluating $g(m,r)$ given
values of $g(m,j)$, $j<r$.  The first theorem provides the means
for evaluating $g(m,r)$ when $r$ belongs to the family of
integers defined by the recurrence relation
$r_{n}=3r_{n-1}-r_{n-2}$ with particular initial conditions.

\begin{theorem}\label{thm:g(m,j)=r(m-1)+1}
Let $m, j$ and $r$ be positive integers, with $m \geq 2$ and $j<r$.  If $r(m-1)\leq g(m,j)\leq
r(m-1)+n$ for $m \geq m_0$, where $r,n$, and $m_0$ are positive
integers, then $$g(m,r)=(3r-j)(m-1)+1$$ for $m \geq
\max\{m_0,n+1\}$.
\end{theorem}
\begin{proof}By hypothesis there exists $\D_j:
[r(m-1)+2,2r(m-1)] \rightarrow [1,j]$ which is an $L(j)$-coloring
for $m \geq m_0$.  For convenience, let
$$\mathcal{I}_i=[(2r+i-1)(m-1)+1,(2r+i)(m-1)]$$ for $i
\in[1,r-j]$. Define the function $\D_r
:[r(m-1)+2,(3r-j)(m-1)]\rightarrow [1,r]$ as follows
\[ \D_r(x) =
   \left\{ \begin{array}{ll}
       \D_j(x), &\mbox{ for }x \in [r(m-1)+2,2r(m-1)] \\
       j+i, &\mbox{ for } x \in \mathcal{I}_i, i
\in[1,r-j]. \end{array} \right. \]

That $\D_r$ is an $L(r)$-coloring follows since $\D_j$ is an
$L(j)$-coloring.  Since for each $c \in \Delta(I_3)$ we have $|\Delta^{-1}(c)\cap (I_2 \cup I_3)|=m-1$, we see that $g(m,r)>(3r-j)(m-1)$ for $m \geq m_0$.

Now, let $\D:[1,(3r-j)(m-1)+1] \rightarrow [1,r]$ be an arbitrary
$r$-coloring and $m \geq \max\{m_0,n+1\}$.  Let $\D(I_3)=C$ and
$k=|C|$. We proceed to show that $\D$ is not an $L(r)$-coloring by
case analysis of $k$.

\textbf{Case 1.}  Suppose $k \leq r-j$.  Since
$|I_3|=(r-j)(m-1)+1$, it follows that there must be some $c \in
[1,r]$ such that $|\D^{-1}(c) \cap I_3| \geq m$, whence $\Delta$ is not an $L(r)$-coloring by Proposition \ref{prop:first}.

\textbf{Case 2.}  Suppose $k > r-j$.  Let $S=\D^{-1}(C)\cap (I_2
\cup I_3)$ and let $U=S \cap I_2$.  Let $t=|\D(I_2)|-|\D(I_2) \cap
C|$, so that $t \leq r-k<j$. Assume without loss of generality that
$\D(I_2) \setminus \{\D(I_2) \cap C\}=[1,t]$.  Furthermore, we may
assume $|S| \leq k(m-1)$, since otherwise some monochromatic
$m$-set $Y$ exists with $Y \subset I_2 \cup I_3$ and $Y \cap I_3
\neq \varnothing$ and we are done.  Hence, since $|I_3|=(r-j)(m-1)+1$, we have
that $|U| =|S|-|I_3| \leq (k-r+j)(m-1)-1$.

Let $\mathcal{P}$ be a partition of $U$ into $p=j-t \geq k-r+j$
sets $\gamma_1,\ldots,\gamma_p$ such that $|\gamma_i|\leq m-1$ for
each $i \in [1,p]$.  Define a coloring $\widehat \D:I_2 \rightarrow
[1,j]$ as follows:

\[\widehat \D(x) =
   \left\{ \begin{array}{ll}
       \D(x), &\mbox{ for }\D(x) \in [1,t] \\
       t+i, &\mbox{ for } x \in \gamma_i, i
\in[1,j-t]. \end{array} \right. \]

Using the notation of Lemma \ref{lem:needsomeofeachcolor} and the fact that $A_c(\D)+B_c(\D) \leq |\D^{-1}(c)|$ when $|\D^{-1}(c)|\leq m-1$, we note
that \begin{equation*}\begin{split}\sum_{c \in [1,j]} A_c(\widehat \D)
+ \min\{B_c(\widehat \D),m-1\} &\leq t(2m-2)+|U| \\ &\leq
(2t+k-r+j)(m-1)-1.\end{split}\end{equation*} Since $g(m,j)-n-1
\leq r(m-1)-1$, Lemma \ref{lem:needsomeofeachcolor} implies that
$\widehat \D$ is not an $L(j)$-coloring if $(2t+k-r+j)(m-1)-1 \leq
j(2m-2)-n-1$.  Using $t \leq r-k<j$, this inequality is easily
verified for $m \geq n+1 \geq 2$. Hence there exist monochromatic
$m$-sets $X,Y \subset I_2$ where $X \prec Y$ and $2(x_m-x_1) \leq
y_m-x_1$. Moreover, $\widehat \D(X) \subseteq [1,t]$ and $\widehat \D(Y) \subseteq
[1,t]$ since $|\widehat \D^{-1}(t+i)|<m$ for each $i \in [1,p]$. Thus,
$X$ and $Y$ are monochromatic in $\D$, and the proof is complete.
\end{proof}

\begin{example}\label{ex:1}Consider the alternate proof that $g(m,2)=5(m-1)+1$ for $m \geq
3$:  note that $g(m,1)$ is trivially $2m = 2(m-1)+2$ for all
positive $m$; by the previous proof, we have $g(m,2) = 5(m-1)+1$
for all $m \geq 3$.

As another example, we have seen in Theorem \ref{thm:g(m,2)} that
$g(m,2)=5(m-1)+1$ for all $m \geq 2$. By the previous theorem,
this implies $g(m,5)=13(m-1)+1$ for $m \geq 2$, which in turn
implies $g(m,13)=34(m-1)+1$ for $m \geq 2$.

Likewise, we have see in Theorem \ref{thm:g(m,4)} that
$g(m,4)=10(m-1)+1$ for all $m \geq 3$.  The previous theorem gives
$g(m,10)=26(m-1)+1$ for $m \geq 3$, which in turn implies
$g(m,26)=68(m-1)+1$ for $m \geq 3$.

More explicitly, Theorems \ref{thm:g(m,2)} and
\ref{thm:g(m,4)} can be used in conjunction with Theorem \ref{thm:g(m,j)=r(m-1)+1} to solve
$g(m,r_n)$, when $r_n$ is in the family of integers generated by
the recurrence relation
\begin{equation}\label{eqn:rec}r_{n}=3r_{n-1}-r_{n-2}\end{equation} with initial
conditions $r_0=2,r_1=5$ from Theorem \ref{thm:g(m,2)} or $r_0=4,r_1=10$ from Theorem \ref{thm:g(m,4)} .

One can solve these recurrence relations in terms of the Fibonacci numbers.  In particular the initial value set $r_0 =2, r_1=5$ gives $r_n = 5f_{2n}-2f_{2n-2}$, where $f_0=0$ and $f_1=1$ are the first two Fibonacci numbers.  Using properties of Fibonacci sequence simplifies this expression to $r_n = f_{2n+3}$.  Of course the recurrence relation with initial conditions $r_0=4$ and $r_1=10$ then has general solution $r_n = 2f_{2n+3}$.
\begin{flushright}$\blacksquare$\end{flushright}\end{example}

Our ultimate goal is to evaluate $g(m,r)$ for as many $r$ as
possible.  Although Theorem \ref{thm:g(m,j)=r(m-1)+1} is an
important step in that direction, it is of no use without the
proper asymptotic value $g(m,r_0)=r_1(m-1)+n$.  We shall need another
result to provide a bound on $g(m,r)$ so that we may apply Theorem
\ref{thm:g(m,j)=r(m-1)+1}.

\begin{theorem}\label{thm:moregoodresults}Let $m, j$ and $r$ be positive integers, with $m \geq 2$ and $j+1<r$.  If $(r-2)(m-1) \leq g(m,j)$ for $m \geq m_0$ and
$g(m,j+1) \leq (r+1)(m-1)+n$ for $m \geq m_1$, where $r,n,m_0$,
and $m_1$ are positive integers, then
$$(3r-j-1)(m-1)<g(m,r) \leq (3r-j-1)(m-1)+n$$ for $m \geq
\max\{m_0,m_1\}$.
\end{theorem}
\begin{proof}
By hypothesis there exists
$\D_j:[(r+1)(m-1)+2,(2r-1)(m-1)]\rightarrow[1,j]$ which is an
$L(j)$-coloring for $m \geq m_0$.  Define
$\D_{j+1}:[r(m-1)+2,2r(m-1)]\rightarrow [1,j+1]$ as follows

\[ \D_{j+1}(x) =
   \left\{ \begin{array}{ll}
       j+1, &\mbox{ for } x \in [r(m-1)+2,(r+1)(m-1)+1]\\ &\mbox{ \ or } x\in [(2r-1)(m-1)+1,2r(m-1)] \\
       \D_j(x), &\mbox{ otherwise. } \end{array} \right. \]

Since $\D_j$ is an $L(j)$-coloring it follows immediately that
$\D_{j+1}$ is an $L(j+1)$-coloring.

As before, let
$$\mathcal{I}_i=[(2r+i-1)(m-1)+1,(2r+i)(m-1)]$$ for $i \in
[1,r-j-1]$.  Define the function $\D_r:I_2 \cup
[2r(m-1)+1,(3r-j-1)(m-1)] \rightarrow [1,r]$ as follows

\[ \D_r(x) =
   \left\{ \begin{array}{ll}
       \D_{j+1}(x), &\mbox{ for } x \in I_2 \\
       j+1+i, &\mbox{ for } x \in \mathcal{I}_i, i
\in[1,r-j-1]. \end{array} \right. \]

From Lemma \ref{lem:colorextension} and its subsequent remark,
$\D_r$ can be extended to an $L(r)$-coloring of
$[1,(3r-j-1)(m-1)]$, and so $g(m,r)>(3r-j-1)(m-1)$.

Let $\D:[1,(3r-j-1)(m-1)+n+1]\rightarrow [1,r]$ be a given
$r$-coloring, and let $m \geq \max\{m_0,m_1\}$.  Let $\D(I_3)=C$
and $k=|C|$. We proceed to show that $\D$ is not an
$L(r)$-coloring by case analysis of $k$.

\textbf{Case 1.}  Suppose $k \leq r-j-1$.  Since
$|I_3|=(r-j-1)(m-1)+n$ where $n \geq 1$ it follows that there must be some $c \in
[1,r]$ such that $|\D^{-1}(c)\cap I_3| \geq m$, whence $\Delta$ is not an $L(r)$-coloring by Proposition \ref{prop:first}.

\textbf{Case 2.}  Suppose $k > r-j-1$.  Let $S=\D^{-1}(C)\cap (I_2
\cup I_3)$ and let $U=S \cap I_2$.  Let $t=|\D(I_2)|-|\D(I_2) \cap
C|$, so that $t \leq r-k<j+1$. Assume without loss of generality that
$\D(I_2) \setminus \{\D(I_2) \cap C\}=[1,t]$.  Furthermore, we may
assume $|S| \leq k(m-1)$, since otherwise some monochromatic
$m$-set $Y$ exists with $Y \subset I_2 \cup I_3$ and $Y \cap I_3
\neq \varnothing$, and we are done.  Since $|I_3|=(r-j-1)(m-1)+n$, we have that
$|U| =|S|-|I_3| \leq (k-r+j+1)(m-1)-n$.

Let $\mathcal{P}$ be a partition of $U$ into $p=j+1-t \geq
k-r+j+1$ sets $\gamma_1,\ldots,\gamma_p$ such that $|\gamma_i|\leq
m-1$ for each $i \in [1,p]$.  Define a coloring $\widehat \D:I_2
\rightarrow [1,j+1]$ as follows:

\[\widehat \D(x) =
   \left\{ \begin{array}{ll}
       \D(x), &\mbox{ for }\D(x) \in [1,t] \\
       t+i, &\mbox{ for } x \in \gamma_i, i
\in[1,j+1-t]. \end{array} \right. \]

Using the notation of Lemma \ref{lem:needsomeofeachcolor} and the fact that $A_c(\D)+B_c(\D) \leq |\D^{-1}(c)|$ when $|\D^{-1}(c)|\leq m-1|$, we have
\begin{equation*}\begin{split}\sum_{c \in [1,j+1]}
A_c(\widehat \D) + \min\{B_c(\widehat \D),m-1\} &\leq t(2m-2)+|U| \\
&\leq (2t+k-r+j+1)(m-1)-n.\end{split}\end{equation*} Since
$g(m,j+1)-m-n \leq r(m-1)-1$, Lemma \ref{lem:needsomeofeachcolor}
implies that $\widehat \D$ is not an $L(j+1)$-coloring if
$(2t+k-r+j)(m-1)-n \leq (j+1)(2m-2)-m-n$. Using $t \leq
r-k<j+1$, this is easily verified for all $m \geq 2$.  Hence there exist
monochromatic $m$-sets $X,Y \subset I_2$ where $X \prec Y$ and
$2(x_m-x_1) \leq y_m-x_1$. Moreover, $\widehat \D(X) \in [1,t]$ and
$\widehat \D(Y) \in [1,t]$ since $|\widehat \D^{-1}(t+i)|<m$ for each $i
\in [1,p]$. Thus, $X$ and $Y$ are monochromatic in $\D$, and the
proof is complete.\end{proof}

\begin{example}\label{ex:2}From Theorems \ref{thm:g(m,2)} and
\ref{thm:g(m,3)} we have that $g(m,2)=5(m-1)+1$ for $m\geq 2$ and
$g(m,3) \leq 8(m-1)+1$ for $m \geq 4$.  We see from Theorem
\ref{thm:moregoodresults} that $g(m,7)=18(m-1)+1$ for $m \geq 4$.  Repeated use of Theorem \ref{thm:g(m,j)=r(m-1)+1} provides another infinite family $\{r_n\}$ for which $g(m,r_n)=r_{n+1}(m-1)+1$. Here the elements $r_n$ satisfy Equation \ref{eqn:rec} with initial conditions $r_0=7,r_1=18$.  This family can also be expressed in terms of the Fibonacci numbers, with $$r_n = 18f_{2n}-7f_{2n-2}.$$

Likewise, from Theorems \ref{thm:g(m,3)} and \ref{thm:g(m,4)} we
have that $g(m,3)>7(m-1)+1$ for $m \geq 4$ and $g(m,4)=10(m-1)+1$
for $m \geq 3$.  Applying Theorem \ref{thm:g(m,j)=r(m-1)+1}, we
have $g(m,9)=23(m-1)+1$ for $m \geq 4$. Again, repeated use of Theorem \ref{thm:moregoodresults} solves $g(m,r_n)=r_{n+1}(m-1)+1$, where here $$r_n = 23f_{2n}-9f_{2n-2}.$$  \begin{flushright}$\blacksquare$\end{flushright}\end{example}

The next result gives a fairly loose bound for $g(m,r)$ given
values of $g(m,j)$, $j<r$.  However, it bounds the function
$g(m,r)$ such that Theorem \ref{thm:moregoodresults} may be
invoked.

\begin{theorem}\label{thm:g(m,r) bounds}Let $m, j$ and $r$ be positive integers, with $m \geq 2$ and $j<r$.  If $(r-1)(m-1)+1 \leq g(m,j) <
r(m-1)$ for $m \geq m_0$, where $r$ and $m_0$ are positive
integers, then
$$(3r-j-1)(m-1)+1 < g(m,r)\leq (3r-j)(m-1)$$ for $m \geq m_0$.
\end{theorem}
\begin{proof}  We start with the lower bound.  By hypothesis there exists
$\D_j:[(r+1)(m-1)+1,2r(m-1)]\rightarrow[1,j]$ which is an
$L(j)$-coloring for $m \geq m_0$.  As before, let
$$\mathcal{I}_i=[(2r+i-1)(m-1)+1,(2r+i)(m-1)]$$ for $i \in
[1,r-j-1]$.  Define the function $\D_r:[r(m-1)+2,(3r-j-1)(m-1)+1]
\rightarrow [1,r]$ as follows

\[ \D_r(x) =
   \left\{ \begin{array}{ll}
       j+1, &\mbox{ if } x \in [r(m-1)+2,(r+1)(m-1)]\\
&\mbox{ \ or } x = (3r-j-1)(m-1)+1 \\
       \D_j(x), &\mbox{ for }x \in [(r+1)(m-1)+1,2r(m-1)] \\
       j+1+i, &\mbox{ for } x \in \mathcal{I}_i, i
\in[1,r-j-1]. \end{array} \right. \]

It is not difficult to see that $\D_r$ is an $L(r)$-coloring on
$I_2$ such that there is no monochromatic $m$-set $Y \subset I_2 \cup I_3$ with
$y_m \in I_3$. Thus, it follows from Proposition \ref{prop:first}
and Lemma \ref{lem:colorextension} that $g(m,r)>(3r-j-1)(m-1)+1$
for every $m \geq m_0$.

To show that $g(m,r)\leq (3r-j)(m-1)$, let $\D:[1,(3r-j)(m-1)]
\rightarrow [1,r]$ be an arbitrary $r$-coloring.  Let $\D(I_3)=C$
and $k=|C|$. We proceed to show that $\D$ is not an
$L(r)$-coloring by case analysis of $k$.

\textbf{Case 1.}  Suppose $k < r-j$.  Since $|I_3|=(r-j)(m-1)$, it
follows that there must be some $c \in [1,r]$ such that
$|\D^{-1}(c)\cap I_3| \geq m$, whence $\Delta$ is not an $L(r)$-coloring by Proposition \ref{prop:first}.

\textbf{Case 2.}  Suppose $k = r-j$.  Since $g(m,j)<r(m-1)$ and
$|I_2|=r(m-1)-1$, if $|\D(I_2)| \leq j$ then $\D$ is not an
$L(j)$-coloring.  Hence $\D(I_2) > j$ so that $\D(I_2) \cap
\D(I_3) \neq \varnothing$, and it follows that there exists some
$z \in \D^{-1}(C)\cap I_2$.  Since $|I_3 \cup
\{z\}|=(r-j)(m-1)+1$, there must be some monochromatic $m$-set $Y$
such that $Y \subset I_2 \cup I_3$ and $y_m \in I_3$.  Applying
Proposition \ref{prop:first} completes the proof.

\textbf{Case 3.}  Suppose $k > r-j$.  Let $S=\D^{-1}(C)\cap (I_2
\cup I_3)$ and let $U=S \cap I_2$.  Let $t=|\D(I_2)|-|\D(I_2) \cap
C|$, so that $t \leq r-k$.  Assume for simplicity that $\D(I_2) \setminus
\{\D(I_2) \cap C\}=[1,t]$.  Furthermore, we may assume $|S| \leq
k(m-1)$, since otherwise some monochromatic $m$-set $Y$ exists
with $Y \subset I_2 \cup I_3$ and $Y \cap I_3 \neq \varnothing$.
Hence, since $|I_3|=(r-j)(m-1)$, we have that $|U| \leq
(k-r+j)(m-1)$.

Let $\mathcal{P}$ be a partition of $U$ into $p=j-t \geq k-r+j$
sets $\gamma_1,\ldots,\gamma_p$ such that $|\gamma_i|\leq m-1$ for
each $i \in [1,p]$.  Define a coloring $\widehat \D:I_2 \rightarrow
[1,j]$ as follows

\[\widehat \D(x) =
   \left\{ \begin{array}{ll}
       \D(x), &\mbox{ for }\D(x) \in [1,t] \\
       t+i, &\mbox{ for } x \in \gamma_i, i
\in[1,j-t]. \end{array} \right. \]

Since $g(m,j)<r(m-1)$ and $|I_2|=r(m-1)-1$, there exist
monochromatic $m$-sets $X,Y \subset I_2$ where $X \prec Y$ and
$2(x_m-x_1) \leq y_m-x_1$.  Moreover, $\widehat \D(X) \in [1,t]$ and
$\widehat \D(Y) \in [1,t]$ since $|\widehat \D^{-1}(t+i)|<m$ for each $i
\in [1,p]$.  Thus, $X$ and $Y$ are monochromatic in $\D$, and the
proof is complete.
\end{proof}

\begin{example}\label{ex:3}By Theorem \ref{thm:g(m,2)} we have
$g(m,2)=5(m-1)+1$ for $m \geq 2$.  Applying Theorem
\ref{thm:g(m,r) bounds} we have $$15(m-1)+1 <
g(m,6)\leq16(m-1)$$ for $m \geq 2$.

Likewise, by Theorem \ref{thm:g(m,3)} we have
$7(m-1)+1 \leq g(m,3)<8(m-1)$ for $m \geq 5$.  From this we see
$$20(m-1)+1 < g(m,8)\leq 21(m-1)$$ for $m \geq 5.$ \begin{flushright} $\blacksquare$ \end{flushright} \end{example}

\section{Conclusion and Conjectures}

In the previous two sections we gave either an exact solution to or a
bound on $g(m,r)$ for all $r \in [2,10]$ and sufficiently large
$m$.  Of course, we could use Theorems \ref{thm:g(m,j)=r(m-1)+1},
\ref{thm:moregoodresults}, and \ref{thm:g(m,r) bounds} to solve or
bound $g(m,r)$ for many $r > 10$.  We conjecture that for each
positive integer $r$ one may find a positive integer $j_r$ such
that one of Theorems \ref{thm:g(m,j)=r(m-1)+1},
\ref{thm:moregoodresults}, or \ref{thm:g(m,r) bounds} may be used
to solve or bound $g(m,r)$.


We have verified by computer the existence of some $j_r$ for each
$r \leq 10^5$.  This program was also used to calculate the
proportions in which exact or bounded results appear in these
first $10^5$ integers, finding that approximately $38.2\%$ of
integers have exact solutions (generated by Theorem
\ref{thm:g(m,j)=r(m-1)+1}), $23.6\%$ are bounded by a constant
(generated by Theorem \ref{thm:moregoodresults}), and the
remaining $38.2\%$ are bounded by a coefficient on $m$ (generated
by Theorem \ref{thm:g(m,r) bounds}).  Furthermore, these
proportions are represented in much smaller samples, perhaps
suggesting that these values are near the asymptotic proportions.

\section*{Acknowledgement}

The author wishes to express his thanks to Professor
A.~Bialostocki for his kind supervision and to D.~Grynkiewicz for
fruitful discussions.  He would also like to thank two anonymous referees for their careful corrections and excellent suggestions.

\end{document}